\documentstyle[verbatim,12pt,leqno]{amsart}

\let\mathcal\cal
\newtheorem{theorem}{Theorem}[section]

\newtheorem{lemma}[theorem]{Lemma}
\newtheorem{corollary}[theorem]{Corollary}

\newtheorem{claim}{Claim}[lemma]
\newtheorem{problem}{Problem}

\theoremstyle{definition}

\newtheorem{definition}[theorem]{Definition}

\theoremstyle{remark}
\newtheorem{remark}{Remark}

\newcommand{\proof}{\begin{pf}}
\newcommand{\Proof}[1]{\begin{pf*}{Proof of #1}}
\newcommand{\eproof}{\end{pf}}
\newcommand{\Eproof}{\end{pf*}}
\newcommand{\sproof}[1]{\begin{pf*}{#1}}
\newcommand{\esproof}{\end{pf*}}

\newcommand{\arablabel}{
          \renewcommand{\labelenumi}{{\rm (\arabic{enumi})}}
          \renewcommand{\theenumi}{{\rm (\arabic{enumi})}}
          \renewcommand{\labelenumii}{{\rm (\arabic{enumii})}}
          \renewcommand{\theenumii}{{\rm (\arabic{enumii})}}
                    }

\newcommand{\alabel}{
          \renewcommand{\labelenumi}{{\rm (\alph{enumi})}}
          \renewcommand{\theenumi}{{\rm (\alph{enumi})}}
          \renewcommand{\labelenumii}{{\rm (\alph{enumii})}}
          \renewcommand{\theenumii}{{\rm (\alph{enumii})}}
                    }

\newcommand{\rlabel}{
          \renewcommand{\labelenumi}{{\rm (\roman{enumi})}}
          \renewcommand{\theenumi}{{\rm (\roman{enumi})}}
          \renewcommand{\labelenumii}{{\rm (\roman{enumii})}}
          \renewcommand{\theenumii}{{\rm (\roman{enumii})}}
                    }

\newcommand{\nolabel}{
         \renewcommand{\labelenumi}{}
         \renewcommand{\theenumi}{} 
         \renewcommand{\labelenumii}{}
         \renewcommand{\theenumii}{} 
                     }

\def\myheads#1;#2;{
\pagestyle{myheadings}
\markboth{{\sc\hfill #1\hfill\protect\makebox[0cm][r]{\rm\today}}}
{{\sc\protect\makebox[0cm][l]{\rm\today}\hfill #2\hfill}}
}

\newcommand{\bcal}{{\mathcal B}}

\newcommand{\dcal}{{\mathcal D}}

\newcommand{\gcal}{{\mathcal G}}

\newcommand{\ical}{{\mathcal I}}

\newcommand{\pcal}{{\mathcal P}}

\newcommand{\ucal}{{\mathcal U}}

\newcommand{\xcal}{{\mathcal X}}

\newcommand{\zcal}{{\mathcal Z}}

\newcommand{\setm}{\setminus}
\newcommand{\empt}{\emptyset}
\newcommand{\subs}{\subset}

\newcommand{\oo}{{{\omega}_1}}
\newcommand{\rest}{\lceil}
\newcommand{\dom}{\operatorname{dom}}
\newcommand{\ran}{\operatorname{ran}}
\def\<{\left\langle}
\def\>{\right\rangle}
\def\cf{\operatorname{cf}}

\def\OO{{\omega}}
\def\oo{\omega_1}
\def\br#1;#2;{\bigl[ {#1} \bigr]^ {#2} }
\def\bc#1;#2;{\bigl( {#1} \bigr)^ {#2} }

\def\ooseq#1;#2;{\< {#1}_{#2}:{#2}<\oo\>}
\def\ooset#1;#2;{\{ {#1}_{#2}:{#2}<\oo\}}
\def\seq#1;#2;#3;{\< {#1}_{#2}:{#2}<#3\>}
\def\set#1;#2;#3;{\{ {#1}_{#2}:{#2}<#3\}}
\def\oseq#1;#2;{\< {#1}_{#2}:{#2}<\OO\>}
\def\oset#1;#2;{\{ {#1}_{#2}:{#2}<\OO\}}
\def\oosequ#1;#2;{\< {#1}^{#2}:{#2}<\oo\>}
\def\oosetu#1;#2;{\{ {#1}^{#2}:{#2}<\oo\}}
\def\sequ#1;#2;#3;{\< {#1}^{#2}:{#2}<#3\>}
\def\setu#1;#2;#3;{\{ {#1}^{#2}:{#2}<#3\}}
\def\osequ#1;#2;{\< {#1}^{#2}:{#2}<\OO\>}
\def\osetu#1;#2;{\{ {#1}^{#2}:{#2}<\OO\}}

\def\force{\raisebox{1.5pt}{\mbox {$\scriptscriptstyle\|$}}
\mbox{$\!\mbox{---}$}}

\newcommand{\Fin}{\operatorname{Fin}}
\newcommand{\fn}{\operatorname{Fn}}

\def\to{\longrightarrow}

\newcommand{\leo}{\operatorname{<_{\text{\rm{On}}}}}

\newcommand{\id}{\operatorname{id}}
\newcommand{\cont}{2^{\omega}}
\def\fin#1;{\br #1;<{\omega};}

\def\w{\operatorname{w}}
\def\nw{\operatorname{nw}}
\def\RR{\operatorname{R}}
\def\D{\operatorname{D}}

\def\rstr{{\rho}^*}
\def\nea{neighbourhood assignment }

\def\anfi#1;{\<A^{#1},n^{#1},f^{#1},g^{#1},k^{#1},d^{#1},e^{#1}\>}
\newcommand{\anfg}{\<A,n,f,g\>}
\def\anfgi#1;{\<A^{#1},n^{#1},f^{#1},g^{#1}\>}
\def\naui#1;{\<A^{#1},\ucal^{#1},f^{#1},g^{#1}\>}
\def\nau{\naui;}
\def\np{n^p}
\def\ap{A^p}
\def\fp{f^p}
\def\up{\ucal^p}
\def\gp{g^p}

\def\aq{A^q}
\def\fq{f^q}
\def\uq{\ucal^q}
\def\gq{g^q}

\def\sbar{\overline{\sigma}}
\def\sstr{{\sigma}^*}

\newcommand{\Pko}{P^{\kappa}}
\newcommand{\Pkoo}{P^{\kappa}_0}
\newcommand{\pcalk}{\pcal^{\kappa}}

\def\ddd#1;#2;{X^#1_#2}
\def\ccc#1;#2;{C^#1_#2}
\newcommand{\dkl}{\ddd {\kappa};{\lambda};}
\newcommand{\cckl}{\ccc {\kappa};{\lambda};}

\newcommand{\lvk}{{\lambda}^{<{\kappa}}}
\newcommand{\vkl}{{}^{<{\kappa}}{\lambda}}
\newcommand{\kl}{{}^{\kappa}{\lambda}}
\newcommand{\lk}{{\lambda}^{\kappa}}
\def\sss#1;{#1^*}

\def\eee#1;#2;{E^{#1}_{#2}}

\newcommand{\pid}{\operatorname{\pi_{\delta}}}
\newcommand{\pin}{\operatorname{\pi_n}}
\newcommand{\px}{\operatorname{\pi_{\xi}}}
\newcommand{\pmm}{\operatorname{\pi_m}}
\newcommand{\swp}{\operatorname{swp}}

\myheads Forcing countable networks $\dots$ 
;
Forcing countable networks $\dots$
;

\title{Forcing countable networks for spaces \\
satisfying $\RR(X^{\omega})={\omega}$}
\thanks{The preparation of this paper was supported by the 
Hungarian National Foundation for Scientific Research grant no. 16391 }

\author{I. Juh\'asz}
\address{Mathematical Institute of the Hungarian Academy of Sciences}
\email{juhasz@@math-inst.hu}
\author{L. Soukup}
\thanks{The second author  was supported by DFG (grant Ko 490/7-1)}
\address{Mathematical Institute of the Hungarian Academy of Sciences}
\email{soukup@@math-inst.hu}
\author{Z. Szentmikl{\'o}ssy}
\address{E{\"o}tv{\"o}s Lor{\'a}nd University, Department of Analysis}

\begin{document}
\maketitle
\begin{center}\large\today\end{center}
\begin{abstract} 	
We show that all finite powers of a Hausdorff space $X$
do not contain uncountable weakly separated subspaces iff
there is a c.c.c poset $P$ such that 
in $V^P$  $X$ is a countable union
of 0-dimensional subspaces of countable weight.
We also show that 
this theorem is sharp in two different senses: (i) we can't get rid of
using generic extensions, (ii) we have to consider all finite powers of $X$.   

\end{abstract}

\section{Introduction}
\label{sc:intr}
We use  standard topological notation and terminology throughout, cf \cite{J}.
The following definitions are less well-known.
\begin{definition}
Given a topological space $\<X,{\tau}\>$ and a 
subspace $Y\subs X$ a function $f$
is called a {\em \nea on Y} iff $f:Y\to {\tau}$ and $y\in f(y)$ for each 
$y\in Y$. 
\end{definition}

\begin{definition}
A space $Y$ is {\em weakly separated} if there is 
a \nea $f$ on $Y$ such that
$$
\forall y\ne z\in Y\ ( y\notin f(z) \lor z\notin f(y) ),
$$
moreover
$$
\RR(X)=\sup\{|Y|:Y\subs X\mbox{ is weakly separated}\}.
$$
\end{definition}

The notion of weakly separated spaces and the cardinal function
$\RR$ were introduced by Tka{\v c}enko in \cite{Tk},
where the following question was also raised:
does $\RR(X^{\omega})={\omega}$ (or even $\RR(X)={\omega}$)
imply that $X$ has a countable network (i.e. $\nw(X)={\omega}$)?
(Note that $\RR(X^{\omega})=\omega$ is equivalent to 
$\RR(X^n)=\omega$ for all $n\in\omega$, moreover 
$\nw(X)=\omega$ implies $\RR(X^\omega)=\nw(X^\omega)=\omega$.)
Several consistent counterexamples  to this were given, e.g. in
\cite{Cie}, \cite{HJ}, \cite{JSSz} and \cite[p. 43]{Stevo2}, but no ZFC counterexample is
known. (In \cite {Stevo2} it is stated that under PFA the implication
is valid, but no proof is given.)
The counterexamples given in \cite{JSSz} and \cite{Stevo2} from CH
are also first countable.

Our main result here says that, at least for  $T_2$ spaces, a weaker version of
Tka{\v c}enko's conjecture is valid, namely $\RR(X^{\omega})={\omega}$
implies that $\nw(X)={\omega}$ holds in a suitable c.c.c and hence
cardinal preserving generic extension!
In fact, in this extension $X$ becomes ${\sigma}$-second countable, i.e. 
$X$ is the union of countably many subspaces of countable weight.

In section \ref{sc:sharp} we show that the main result is sharp in
different senses. Firstly, we  force, for every natural number $n$,
a 0-dimensional, first countable space $X$ such that 
$\RR(X^n)={\omega}$, but $\nw(X)>{\omega}$ in any cardinal preserving
extension of the ground model. 
Secondly, we construct  in ZFC a 0-dimensional 
$T_2$ space $X$ such that $\chi(X)=\nw(X)=\omega$ but $X$ is not
$\sigma$-second countable. 

It easily follows from the proof of our main result 
that if $MA_{\kappa}$ holds and 
$X$ is a Hausdorff  space with $|X|+\w(X)\le {\kappa}$, then
$\RR(X^{\omega})={\omega}$ if and only if $X$ is ${\sigma}$-second
countable. We also prove that $MA(Cohen)$ is not enough
to yield this equivalence. To do this we use a result of
Shelah (the proof of which  presented here with his kind permission)
saying that in any generic extension by Cohen reals
 the ideal of the first category subspaces of a space 
from the ground model  is generated by the 
subspaces of first category from the ground model. 

\section{The main result}\label{sc:equiv}
\begin{theorem}\label{tm:main}
Given a Hausdorff topological space $X$ the following are equivalent:
\begin{enumerate}\arablabel
\item $\RR(X^{\omega})={\omega}$,
\item there is a c.c.c poset $P$ such that $V^P\models$
``$\nw(X)={\omega}$'',
\item there is a c.c.c poset $P$ such that \newline
$V^P\models$ ``{$X$ is a countable union of 
0-dimensional subspaces of countable weight.}''
\end{enumerate}
\end{theorem}

\proof
Since the implications (3) $\Rightarrow$ (2)  $\Rightarrow$ (1) are clear, 
it remains to  prove only  that (1) implies (3).
So assume that $X$ is Hausdorff with  $\RR(X^{\omega})={\omega}$,
and fix a base $\bcal$ of $X$ and   a well-ordering $\prec$ on 
$X\cup\bcal$.

The space $X$ contains at most countably many isolated points
by $\RR(X)={\omega}$, so it is enough to force an appropriate partition 
of $X'$, 
the set of non-isolated points of $X$. 
We say that a 4-tuple $\nau$ is in $P$ provided  (i)--(v) below hold:
\begin{enumerate}\rlabel
\item  $A\in\br X';<{\omega};$ and 
$\ucal\in \br \bcal;<{\omega};$,
\item $f$ and $g$ are functions,
\item $f:A\to {\omega}$, $g:\ucal\to {\omega}\times {\omega}$,
\item  if $g(U)=\<n,i\>$ and $f(x)=n$ 
then $x\notin(\overline{U}\setm U)$
whenever $x\in A$ and $U\in \ucal$, 
\item if $g(U)=g(V)=\<n,i\>$ and $f(x)=n$ 
then $x\in U$ iff $x\in V$ whenever $x \in A$ and $U,V\in \ucal$.
\end{enumerate}

Our idea is that  $f$ will guess the partition of $X'$ into countably many
pieces, $\{F_n:n<{\omega}\}$;
if $g(U)=\<n,i\>$ then  $U\cap F_n$ will be clopen in the subspace
$F_n$, and
$g(U)=g(V)=\<n,i\>$ implies  $U\cap F_n=V\cap F_n$.
Consequently, each $F_n$ will have a countable clopen base.

For $p\in P$ we write $p=\naui p;$.
If $p,q\in  P$ we set $p\le q$ iff 
$\fp\supset \fq$ and  $\gp\supset \gq$.

Two conditions $p$ and $q$ from $P$ are called {\em twins} provided 
 $|\ap|=|\aq|$,
$|\up|=|\uq|$ and denoting by ${\eta}$ and by ${\rho}$
the unique $\prec$-preserving bijections between $\ap$ and $\aq$, 
and between
$\up$ and $\uq$, respectively, we have 
\begin{enumerate}\arablabel
\item ${\eta}\rest \ap\cap \aq=\id$,
 ${\rho}\rest \up\cap \uq=\id$,
\item $\forall x\in \ap$   $\fp(x)=\fq({\eta}(x))$.
\item   $\forall U\in\up$ $\gp(U)=\gq({\rho}(U))$.
\item $\forall x\in \ap$ $\forall U\in\up$
( $x\in U$ iff ${\eta}(x)\in {\rho}(U)$, 
and $x\in \overline{U}$ iff ${\eta}(x)\in \overline{{\rho}(U)}$ ). 
\end{enumerate}

\begin{lemma}\label{lm:ccc}
$\pcal=\<P,\le\>$ satisfies c.c.c.
\end{lemma}

\Proof{lemma \ref{lm:ccc}}
Let $\ooset p;{\alpha};\subs P$. Write 
$p_{\alpha}=\naui {\alpha};$.
Using standard $\Delta$-system and counting arguments  we can assume that
these conditions  are pairwise twins. 
Let $k=|A^{\alpha}|$ and  
$\{a_{{\alpha},i}:i<k\}$ be the $\prec$-increasing enumeration
of $A^{\alpha}$.  For each ${\alpha}<\oo$ and $i<k$ put
\begin{gather}
\notag\ucal^0_{{\alpha},i}=\{U\in\ucal^{\alpha}:
\exists j\ g^{\alpha}(U)=\<(f^{\alpha}(a_{{\alpha},i}),j\>
\land a_{{\alpha},i}\in U\},\\
\notag\ucal^1_{{\alpha},i}=\{U\in\ucal^{\alpha}:
\exists j\ g^{\alpha}(U)=\<(f^{\alpha}(a_{{\alpha},i}),j\>
\land a_{{\alpha},i}\notin U\},\\
\intertext{and finally}\notag
V_{{\alpha},i}=\bigcap\{U:U\in\ucal^0_{{\alpha},i}\}\cap
\bigcap\{X\setm \overline{U}:U\in\ucal^1_{{\alpha},i}\}. 
\end{gather}
By (iv) we have $a_{{\alpha},i}\notin \overline{U}$
for $U\in \ucal^1_{{\alpha},i}$, so 
 $a_{{\alpha},i}\in V_{{\alpha},i}$.
Since $\RR(X^k)={\omega}$ there are ${\alpha}<{\beta}<\oo$ such that 
\begin{equation}
\tag{\dag}
\text{$a_{{\alpha},i}\in V_{{\beta},i}$ and
$a_{{\beta},i}\in V_{{\alpha},i}$ for each $i<k$.}
\end{equation}
We claim that $p_{\alpha}$ and $p_{\beta}$ are compatible in $P$.
Let ${\eta}$ and ${\rho}$ be the functions witnessing that
$p_{\alpha}$ and $p_{\beta}$ are twins. 
Put  $A=A^{\alpha}\cup A^{\beta}$, 
$\ucal=\ucal^{\alpha}\cup \ucal^{\beta}$, $f=f^{\alpha}\cup f^{\beta}$, $g=g^{\alpha}\cup g^{\beta}$
and $p=\nau$. Since $p_{\alpha}$ and $p_{\beta}$  are twins,  
$p$ satisfies (i)--(iii) and $p\le p_{\alpha}, p_{\beta}$.
So all we have to do is to show that  $p$ satisfies (iv) and (v).

\renewcommand{\theclaim}{}
\begin{claim}
If  $U\in\ucal^0_{{\alpha},i}\cup\ucal^1_{{\alpha},i}$, then
$a_{{\alpha},i}\in U$ iff $a_{{\alpha},i}\in {\rho}(U)$. 
\end{claim}

\Proof{the claim}
We know   $a_{{\alpha},i}\in U$ iff $a_{{\beta},i}\in {\rho}(U)$.
Thus  $a_{{\alpha},i}\in U$ implies that ${\rho}(U)\in \ucal^0_{{\beta},i}$
and so $a_{{\alpha},i}\in V_{{\beta},i}\subs {\rho}(U)$.
On the other hand, if $a_{{\alpha},i}\notin U$, then 
$a_{{\beta},i}\notin \rho(U)$, hence 
${\rho}(U)\in \ucal^1_{{\beta},i}$,
and so $a_{{\alpha},i}\in V_{{\beta},i}\subs X\setm\overline{{\rho}(U)}$,
 i.e.,  $a_{{\alpha},i}\notin \overline{{\rho}(U)}$.
\Eproof

Now we check (iv).
Assume   $a_{{\alpha},i}\in A^{\alpha}$,
 $V\in \ucal^{\beta}$ with  
 $g^{\beta}(V)=\<f^{\alpha}(a_{{\alpha},i}),l\>$
and $a_{{\alpha},i}\in\overline{V}$.
We have to  show that $a_{{\alpha},i}\in V$.
Since  $a_{{\alpha},i}\in
V_{{\beta},i}$ by $(\dag)$, we have
  $V\notin \ucal^1_{{\beta},i}$. 
But $f^{\alpha}(a_{\alpha,i})=f^{\beta}(a_{\beta,i})$ since $p_\alpha$ and $p_{\beta}$ are twins, hence
$g^{\beta}(V)=\<f^{\beta}(a_{{\beta},i}),l\>$
implies $V\in\ucal^0_{{\beta},i}\cup \ucal^1_{{\beta},i}$, so
$V\in \ucal^0_{{\beta},i}$. 
Hence   $a_{{\alpha},i}\in V_{\beta,i}\subs V$ by $(\dag)$, which
was to be proved.

Finally we check (v).
Assume that  $a_{{\alpha},i}\in A^{\alpha}$ and 
$U,V\in\ucal^{\alpha}\cup\ucal^{\beta}$ are such that 
$g(U)=g(V)=\<f(a_{{\alpha},i}),l\>$.
Define the function $\rstr:\ucal\to \ucal^{\alpha}$ by the formula
$\rstr=\id\rest{\ucal^{\alpha}}\cup ({\rho})^{-1}$
Then $a_{{\alpha},i}\in U$ iff $a_{{\alpha},i}\in \rstr(U)$ and
$a_{{\alpha},i}\in V$ iff $a_{{\alpha},i}\in \rstr(V)$ by the previous
claim.

But $g(\rstr(U))=g(U)=g(V)=g(\rstr(V))=\<f(a_{{\alpha},i}),l\>$, so 
$a_{{\alpha},i}\in\rstr(U)$ iff $a_{{\alpha},i}\in\rstr(V)$
for $p_{\alpha}$ satisfies (v). Thus 
$a_{{\alpha},i}\in U$ iff 
$a_{{\alpha},i}\in\rstr(U)$ iff
$a_{{\alpha},i}\in\rstr(V)$ iff
$a_{{\alpha},i}\in V$, which proves (v).
\Eproof 

Now let $\gcal$ be a $P$-generic filter and
let $F=\bigcup\{f^p:p\in\gcal\}$ and
$G=\bigcup\{g^p:p\in\gcal\}$.
For $n<{\omega}$ let $F_n=F^{-1}\{n\}$.

\begin{lemma}
$\dom(F)=X'$ and  $\dom(G)=\bcal$.
\end{lemma}

\proof
Let $p=\nau\in P$, $x\in X'\setm A$ and $U\in \bcal\setm \ucal$.
Let $A^*=A\cup\{x\}$, $\ucal^*=\ucal\cup \{U\}$ and
$n=\max\ran(f)+1$. 
Let $\dom(f^*)=A^*$, $f^*\supset f$ and $f^*(x)=n$, $\dom(g^*)=\ucal^*$,
$g^*\supset g$ and $g^*(U)=\<n+1,0\>$.  
Then it is easy to check  $p^*=\naui *;\in P$ and obviously 
$p^*\le p$. So the lemma holds because
a generic filter intersects every dense set.
\eproof 

For $m\in {\omega}$ let $$
\bcal_m=\{U\cap F_m:U\in\bcal\text{ and } G(U)=\<m,i\> 
\text{for some } i\in {\omega} \}.
$$

\begin{lemma}
$\bcal_m$ is a countable, clopen base of the subspace $F_m$ of $X'$.
\end{lemma}

\proof
If $U\in\bcal$,  $G(U)=\<m,i\>$ then $U\cap F_m$ is clopen
in $F_m$ by (iv). 
If $U,V\in\bcal$, 
$G(U)=G(V)=\<m,i\>$ then $U\cap F_m=V\cap F_m$ by (v).
So $\bcal_m$ is countable. 

Finally we show that it is a base of $F_m$. So fix $x\in F_m$ and  $V\in\bcal$ 
with $x\in V$. Let $p=\nau\in P$ such that $f(x)=m$. 
 Since $X$ is Hausdorff and  $x$ is non-isolated in $X$,  
we can choose $U\in\bcal\setm \ucal $ such that $x\in U\subs V$
and $\overline{U}\cap A=\{x\}$.

Let  $\ucal^*=\ucal\cup \{U\}$. 
Define the function  
$g^*:\ucal^*\to {\omega}\times {\omega}$ such that $g^*\supset g$ and
 $g^*(U)=\<m,k\>$ where $k=\min\{l:\ran g\subs l\times l\}$.
Then $p^*=\<A,\ucal^*,f,g^*\>$ is an extension of $p$ in $\pcal$
and $p^*\force x\in U\cap F_m\subs V\cap F_m\land U\cap F_m\in\bcal_m$.
Consequently, if $p\force$ ``$x\in F_m\cap V$'' then 
we also have $p\force$ ``$\exists U\in\bcal_m (x\in U\subs F_m\cap V)$'',
which completes the proof.
\eproof
Thus theorem  \ref{tm:main}  is proved.
\eproof

It is easy to check that the above proof needs the
genericity of $\gcal$ over $|X|+|\bcal|$ many
dense sets only, and this immediately yields the following result.

\begin{corollary}\label{cor:ma}
If $MA_{\kappa}$ holds then  for a Hausdorff space $X$ with
$|X|+\w(X)\le {\kappa}$ the following are equivalent:
\begin{enumerate}\arablabel
\item $\RR(X^{\omega})={\omega}$,
\item $\nw(X)={\omega}$,
\item $X$ is the union of  countably many
0-dimensional subspaces of countable weigh,
\item $X$ is the union of  countably many
separable metrizable  subspaces.
\end{enumerate}
\end{corollary}

\section{Sharpness of the main result}\label{sc:sharp}
Our aim in this section is to examine how sharp the above main result is.
The co-finite topology on any uncountable set $X$ clearly satisfies 
$\RR(X^\omega)=\omega$, while, in any extension,
$\nw(X)=|X|$. This show that in the proof of 
$(1)\to (2)$ of \ref{tm:main} the Hausdorffness of $X$ cannot be 
replaced by $T_1$. 
 
The next result in this section  implies that, at least in ZFC, 
the exponent ${\omega}$ in proving $(1)\to (2)$ in 
theorem \ref{tm:main} can not be lowered.

\begin{theorem}
\label{tm:rn}
For each uncountable cardinal ${\kappa}$ and natural number $m$
there is a c.c.c poset $\pcal$
of cardinality ${\kappa}$ such that in $V^{\pcal}$ 
there is a   0-dimensional
first countable topological space $X=\<{\kappa},{\tau}\>$ such that
$\RR(X^m)={\omega}$ but $\RR(X^{m+1})={\kappa}$, hence 
$\nw(X)=\kappa$ in any cardinal preserving extension.
\end{theorem}

In \cite[theorem 3.5]{JSSz} we constructed a c.c.c poset $\<\Pko,\le\>$
which adds to the 
ground model a 0-dimensional, first countable topology $\tau$ on $\kappa$ 
such that $\RR(X^{\omega})={\omega}$ and $\w(X)={\kappa}$
for $X=\<\kappa,\tau\>$. 
The conditions in $\Pko$ are  finite approximations 
of the space $X$ and the property $\RR(X^{\omega})={\omega}$ is guaranteed
by some $\Delta$-system and amalgamation arguments.
Here we will use a subset $P$ of $\Pko$ with the inherited
order. 
To ensure $\RR(X^{m+1})={\kappa}$  we thin out $\Pko$ in the following way.
We fix a family $\dcal=\{d_{\alpha}:{\alpha}<{\kappa}\}$
 of pairwise disjoint elements of ${\kappa}^{m+1}$ with the intention
to make $\dcal$ discrete in $X^{m+1}$.  
A condition $p\in \Pko$ is   put into $P$ 
if and only if every neighbourhood given by $p$ witnesses that 
$\dcal$ is discrete. The main step of the proof is to show
that  $P$ is large enough to allow 
the $\Delta$-system and amalgamation arguments to work in 
showing $\RR(X^m)={\omega}$.

\Proof{theorem \ref{tm:rn}}
First  we  recall some definitions and lemmas from  the proof of  
\cite[theorem 3.5]{JSSz}. 
A quadruple $\anfg$ is said to be in $\Pkoo$ provided
$(a)$--$(b)$ below hold:
\begin{enumerate}\alabel
\item $A\in\br {\kappa};<{\omega};$, $n\in {\omega}$,
$f$ and  $g$ are functions,
\item $f:A\times A\times n\to 2$,
$g:A\times n\times A\times n\to 3$,
\end{enumerate}

For $p\in \Pkoo$ we write $p=\anfgi p;$.
If $p, q\in \Pkoo$ we set 
$p\le q$ iff   $f^p\supseteq f^q$ and $g^p\supseteq g^q$. 
If $p\in \Pkoo$, ${\alpha}\in \ap$, $i<\np$ we defined
$
U({\alpha},i)=U^p({\alpha},i)=\{{\beta}\in \ap:\fp({\beta},{\alpha},i)=1\}
$.

A quadruple $\anfg\in \Pkoo$ is put  in $\Pko$ iff (i)--(ii) below are also 
satisfied:
\begin{enumerate}\rlabel
\item $\forall {\alpha}\in A$ $\forall i<j<n$  
${\alpha}\in U({\alpha},j)\subs U({\alpha},i)$,
\item $\forall {\alpha}\ne {\beta}\in A$ $\forall i,j<n$
\newline
\begin{tabular}{rl}
$g({\alpha},i,{\beta},j)=0$&if and only if $U({\alpha},i)\subs U({\beta},j)$,
\\
$g({\alpha},i,{\beta},j)=1$&if and only if 
$U({\alpha},i)\cap U({\beta},j)=\empt$,
\\
$g({\alpha},i,{\beta},j)=2$&
if ${\alpha}\in U({\beta},j)$ and ${\beta}\in U({\alpha},i)$.
\end{tabular}
\end{enumerate}

\begin{definition}(\cite[definition 3.6]{JSSz})
Assume that $p_i=\anfgi i;\in \Pkoo$ for $i\in 2$. We say that $p_0$ and $p_1$
are {\em twins} iff  $n_0=n_1$, $|A_0|=|A_1|$ and taking $n=n_0$ and
denoting by ${\sigma}$ the unique $\leo$-preserving bijection between $A_0$ and
$A_1$ we have
\begin{enumerate}\arablabel
\item ${\sigma}\rest {A_0\cap A_1}=\id_{A_0\cap A_1}$.
\item ${\sigma}$ is an isomorphism between $p_0$ and $p_1$, i.e.
$\forall {\alpha},{\beta}\in A_0$, $\forall i,j<n$
\begin{enumerate}\nolabel
\item $f_0({\alpha},{\beta},i)=f_1({\sigma}({\alpha}),{\sigma}({\beta}),i)$,
\item $g_0({\alpha},i,{\beta},j)=
g_1({\sigma}({\alpha}),i,{\sigma}({\beta}),j)$,
\end{enumerate}
\end{enumerate}
We say that ${\sigma}$ is the {\em twin function} of $p_0$ and $p_1$.
Define the {\em smashing function} $\sbar$ of $p_0$ and $p_1$ as follows:
 $\sbar={\sigma}\cup \id_{A_1}$.  
The function $\sstr$ defined by the formula 
$\sstr={\sigma}\cup {\sigma}^{-1}\rest {A_1}$ is called the
{\em exchange function} of $p_0$ and $p_1$.
\end{definition}

\begin{definition}(\cite[definition 3.7]{JSSz})
\label{def:eps-amalg}
Assume that $p_0$ and $p_1$ are twins and 
${\varepsilon}:A^{p_1}\setm A^{p_0}\to 2$. A common extension 
$q\in\Pko$ of $p_0$ and $p_1$  is called an
{\em ${\varepsilon}$-amalgamation} of the twins $p_0$ and $p_1$ provided
$$\forall {\alpha}\in A^{p_0}\triangle A^{p_1}
\ f^q({\alpha},\sstr({\alpha}),i)={\varepsilon}(\sbar({\alpha})).
$$ 
\end{definition}


\begin{lemma} (\cite[lemma 3.8]{JSSz})
\label{lm:old}
If $p_0$,  $p_1\in\pcalk$ are twins and 
${\varepsilon}:A^{p_1}\setm A^{p_0}\to 2$, then $p_0$
and $p_1$ have an ${\varepsilon}$-amalgamation in $\Pko$.
\end{lemma}

In \cite{JSSz} we used the poset 
$\pcalk=\<\Pko,\le\>$.
Here we will apply a subset $P$ of $\Pko$.
To define it let  
$\{d_{\alpha}:{\alpha}<{\kappa}\}$
be a family of pairwise disjoint elements of $\br {\kappa};m+1;$
such that 
${\kappa}\setm \bigcup\{d_{\alpha}:{\alpha}<{\kappa}\}$ is still
infinite.  Write $d_{\alpha}=\{d_{{\alpha},i}:i\le m\}$.

\begin{definition}\label{df:pcal}
A condition $p=\anfg\in \Pko$ is in $P$ iff it satisfies (1) and (2) below:
\begin{enumerate}\arablabel
\item $d_{\alpha}\subs A$ or $d_{\alpha}\cap A=\empt$ for each 
${\alpha}<{\kappa}.$
\item if ${\alpha}<{\beta}<{\kappa}$, 
$d_{\alpha}\cup d_{\beta}\subs A$,
then there is an $i=i_{{\alpha},{\beta}}\le m$
such that 
$$
d_{{\alpha},i}\notin U^p(d_{{\beta},i},0)\text{ and }
d_{{\beta},i}\notin U^p(d_{{\alpha},i},0).
$$ 
\end{enumerate}
Let $\pcal=\<P,\le\>$.
\end{definition}

We define $X$ as expected.
Let $\gcal$ be a $\pcal$-generic filter and
let $F=\bigcup\{f^p:p\in\gcal\}$.
For each ${\alpha}<{\kappa}$ and $n\in{\omega}$ let
$V({\alpha},i)=\{{\beta}<{\kappa}:F({\beta},{\alpha},i)=1\}$.
Put $\bcal_{\alpha}=\{V({\alpha},i):i<{\kappa}\}$ and
$\bcal=\bigcup\{\bcal_{\alpha}:{\alpha}<{\kappa}\}$.
We choose $\bcal$ as the base of $X=\<{\kappa},{\tau}\>$.
By standard density arguments we can see that $X$
is first countable and  0-dimensional.

It is easy to see that
$\RR(X^{m+1})={\kappa}$, in fact  $\operatorname{s}(X^{m+1})={\kappa}$.
Indeed, by \ref{df:pcal}(2), $\{d_{\alpha}:{\alpha}<{\kappa}\}$
is discrete in $X^{m+1}$, as  witnessed by the open neighborhoods 
$V(d_{{\alpha},0},0)\times V(d_{{\alpha},1},0)\dots  
\times V(d_{{\alpha},m},0)$.

Finally we need to show that  $\pcal$ satisfies c.c.c
and $V^{\pcal}\models\RR(X^m)={\omega}$.
Clearly, both of these statement follow from the next lemma.

\begin{lemma}\label{lm:ccc+rm}
If $\{p_{\gamma}:{\gamma}<\oo\}\subs \pcal$,
$\{c_{\gamma}:{\gamma}<{\omega}_1\}\subs{\kappa}^m$ 
and $j_0,\dots, j_{m-1}$ are natural numbers, then there 
are ordinals $\alpha<\beta<\oo$ and a condition $p\in P$ such that  
\begin{enumerate}
\item[(+)] $p\le p_{\alpha},p_{\beta}$ and 
$
p\force c_{\alpha}\in 
\prod_{i<m}V(c_{{\beta}}(i),j_i) \land
c_{{\beta}}\in \prod_{i<m}V(c_{{\alpha}}(i),j_i). 
$
\end{enumerate}
\end{lemma}
   
\proof
We can assume that $c_{\gamma}\subs A^{p_{\gamma}}$ holds for 
each ${\gamma}<\oo$. 

Pick ${\alpha}$ and ${\beta}$ such that $p_{\alpha}$ and $p_{\beta}$
are twins, and denoting by ${\rho}$ their twin function we have
 ${\rho}''c_{\alpha}=c_{\beta}$ and  
$\{{\rho}''d_{\xi}:d_{\xi}\subs A^{p_{\alpha}}\}=
\{d_{\zeta}:d_{\zeta}\subs A^{p_{\beta}}\}$.

Let $x=\{{\xi}<{\kappa}:d_{\xi}\subs A^{p_{\alpha}}\}$
and $y=\{{\xi}<{\kappa}:d_{\xi}\subs A^{p_{\beta}}\}$.
Define the function ${\varepsilon}:A^{p_{\alpha}}\setm A^{p_{\beta}}\to 2$ 
by the stipulations ${\varepsilon}({\nu})=1$ iff ${\nu}\in c_{\alpha}$.
By lemma \ref{lm:old} $p_{\alpha}$ and 
$p_{\beta}$ has an ${\varepsilon}$-amalgamation $p$ in $\Pko$. 
Since $C={\kappa}\setm \bigcup\{d_{\xi}:{\xi}<{\kappa}\}$ is infinite,
we can assume that 
$A^p\setm (A^{p_{\alpha}}\cup A^{p_{\beta}})\subs C$.

First we show that $p\in P$. 
Observe that for any ${\xi}<{\kappa}$ we have 
$d_{\xi}\cap A^p\ne\empt$ if and only if 
$d_{\xi}\cap (A^{p_{\alpha}}\cup A^{p_{\beta}})\ne\empt$ if and only if 
$(d_{\xi}\subs A^{p_{\alpha}}$ $\lor$ $d_{\xi}\subs A^{p_{\beta}})$ 
by $A^p\setm (A^{p_{\alpha}}\cup A^{p_{\beta}})\subs C$.
Thus \ref{df:pcal}(1) holds. To check \ref{df:pcal}(2)
assume that ${\xi}\ne{\zeta}<{\kappa}$ and $d_{\xi}\cup d_{\zeta}\subs A^p$.
Then $d_{\xi}\cup d_{\zeta}\subs A^{p_{\alpha}}\cup A^{p_{\beta}}$
and since $p_{\alpha}$ and $p_{\beta}$ are in $P$ we can assume that
$d_{\xi}\subs A^{p_{\alpha}}$ and $d_{\zeta}\subs A^{p_{\beta}}$.
If $d_{\xi}\cap A^{p_{\beta}}\ne\empt$ then $d_{\xi}\subs A^{p_{\beta}}$
as well. Therefore $d_{\xi}\cup d_{\zeta}\subs A^{p_{\beta}}$, and so
\ref{df:pcal}(2) holds for ${\xi}$ and ${\zeta}$ because $p_{\beta}\in P$.
Thus we can assume that $d_{\xi}\subs A^{p_{\alpha}}\setm A^{p_{\beta}}$,
and similarly that $d_{\zeta}\subs A^{p_{\beta}}\setm A^{p_{\alpha}}$.
If ${\rho}''d_{\xi}\ne d_{\zeta}$, then 
${\rho}''d_{\xi}=d_{\mu}$ for some ${\mu}\in y\setm\{{\zeta}\}$.
Since $p_{\beta}\in P$,   there is $i\le m$ such that 
$d_{{\mu},i}\notin U^{p_{\beta}}(d_{{\zeta},i},0)$ and
$d_{{\zeta},i}\notin U^{p_{\beta}}(d_{{\mu},i},0)$. 
So, by the definition of ${\varepsilon}$-amalgamation,
$d_{{\xi},i}\notin U^p(d_{{\zeta},i},0)$ and
$d_{{\zeta},i}\notin U^p(d_{{\xi},i},0)$. 
On the other hand, if ${\rho}''d_{\xi}= d_{\zeta}$, 
then there is $i\le m$ such that 
${\varepsilon}(d_{{\xi},i})=0$, for 
$|d_{\xi}|=m+1>m=|{\varepsilon}^{-1}\{1\}|$.
So, by the definition of ${\varepsilon}$-amalgamation,
$d_{{\xi},i}\notin U^p(d_{{\zeta},i},0)$ and
$d_{{\zeta},i}\notin U^p(d_{{\xi},i},0)$. Thus
$p\in P$.

Finally we show that 
$p\force$ ``$c_{{\alpha}}\in 
\prod_{i<m}V(c_{{\beta}}(i),j_i) \land
c_{\beta}\in  \prod_{i<m}V(c_{{\alpha}}(i),j_i)
$''.
Indeed,   $c_{{\beta}}(i)\in U^p(c_{{\alpha}}(i),j_i)$
and $c_{{\alpha}}(i)\in U^p(c_{{\beta}}(i),j_i)$ for each $i<m$
because ${\rho}(c_{{\alpha},i})=c_{{\beta},i}$ and either 
$c_{{\alpha},i}=c_{{\beta},i}$ or
${\varepsilon}(c_{{\alpha},i})=1$.
\eproof
Theorem \ref{tm:rn} is proved.
\Eproof

Next we show  that the use  of forcing in the implications
$(1)\to (3)$ and $(2)\to (3)$ from
theorem \ref{tm:main} is essential because in \ref{tm:need-forcing} 
we shall produce a ZFC example of a 0-dimensional, 
first countable space $X$ that satisfies 
$\nw(X)=\omega$ (hence $\RR(X^\omega)=\omega$) but still
$X$ is not $\sigma$-second countable. 
To achieve this we need the following lemma.
If $X=\<X,{\tau}\>$ is a topological space, $\D(X)$ denotes 
the discrete topology on  $X$. 
If $A$ and $B$ are sets, let $\Fin(A,B)$ be the family of functions
mapping a finite subset of $A$ into $B$.

\begin{lemma}\label{lm:somewhere-dense}
If  $Z\subs X^{\omega}$ and $Z$ is  somewhere dense 
in the  space $\D(X)^{\omega}$, then $\w(Z)=\w(X)$.
\end{lemma}

\proof
Fix a natural number $n\in {\omega}$ and a function $f:n\to X$
such that  $Z$  is dense in the basic open set  
$U_f=\{g\in X^{\omega}:f\subs g\}$ of $\D(X)^{\omega}$.
This means that 
\begin{equation}
\tag{\dag}\forall f'\in\Fin({\omega}\setm n,X)
\exists g\in Z\ f\cup f'\subs g.
\end{equation}

From now on we forget about the $\D(X)^{\omega}$ topology, we will
use only $(\dag)$.
Without loss of generality we can assume that $Z\subs U_f$.
Let $\zcal$ be a base of $Z$ in the subspace topology of $X^{\omega}$.

Let $\pmm:X^\omega\to X$ be the projection to the 
$m^{\rm th}$ factor, i.e. $\pmm(g)=g(m)$.
Set 
$\xcal=\{\pin(U):U\in\zcal\}$.
Since $\w(Z)\le \w(X)$, it is enough to show that $\xcal$ is a base
of $X$.
\renewcommand{\theclaim}{}
\begin{claim}
If $U$ is open in $Z$ then $\pin(U)$ is 
open in $X$.
\end{claim}

\Proof{the claim }
Let $x\in\pin(U)$. We need to show that 
$\pin(U)$ contains a neighbourhood of $x$.  
Pick $g\in U$ with $x=g(n)$. 
Then, by the definition of the product topology on $X^{\omega}$,
there is a function ${\sigma}$ which maps a finite subset of
${\omega}\setm n$ into the family of non-empty open subsets of
$X$ such that 
\begin{equation}\tag{$\star$}
g\in Z\cap \bigcap_{m\in\dom({\sigma})}\pmm^{-1}{\sigma}(m)\subs U.  
\end{equation}
We can assume that $n\in\dom({\sigma})$.
Let $g'=g\rest(\dom {\sigma}\setm \{n\})$.
By $(\dag)$, for each $x'\in {\sigma}(n)$
there is $h_{x'}\in Z$ such that 
\begin{equation}
\tag{$\star\star$} f\cup\{\<n,x'\>\}\cup f'\subs h_{x'}.
\end{equation}
Now $(\star)$ implies $h_{x'}\in U$ and so $x'\in\pin(U)$.
Thus $x\in {\sigma}(n)\subs \pin(U)$, which was 
to be proved.
\Eproof

To show that $\xcal$ is a base let $x\in V\subs X$, $V$ open.
By $(\dag)$ we can find  a point $g\in Z$ with $f\subs g$ and $g(n)=x$.
The family $\zcal$ is a base of $Z$ in the subspace topology of 
$X^{\omega}$, so there is $U\in \zcal$ such that
$g\in U\subs Z\cap \pin^{-1}V$. Thus $x\in \pin(U)\subs V$ and 
$\pin(U)$ is open by the previous claim.

Thus $\xcal$ is a base of $X$, and so 
$\w(X)\le|\xcal|\le|\zcal|$. Since 
$\w(Z)\le {\omega}\w(X)=\w(X)$, we are done.  
\eproof

After this preparation we can give the  ZFC example promised above.

\begin{theorem}\label{tm:need-forcing}
There is 0-dimensional Hausdorff space  $Y$ such that 
${\chi}(Y)\nw(Y)={\omega}$,
but $Y$ is not ${\sigma}$-second countable.
\end{theorem}

\proof
By \cite[theorem 3.1]{JSSz} there is a 0-dimensional Hausdorff space 
$X$ of size $\cont$ 
such that ${\chi}(X)\nw(X)={\omega}$, but
$\w(X)=\cont$. We show that $Y=X^{\omega}$ is as required.
Clearly ${\chi}(Y)=\chi(X)={\omega}$ and $\nw(Y)=\nw(X)={\omega}$.

Assume that $Y=\bigcup_{k<{\omega}}Z_k$. 
The Baire category theorem implies that some $Z_k$ is somewhere dense 
in   $\D(X)^{\omega}$.
Then $\w(Z_k)=\w(X)=\cont$ by lemma \ref{lm:somewhere-dense}.
\eproof


By corollary \ref{cor:ma}, if Martin's Axiom holds, 
then every Hausdorff space $X$ of  size and weight $<2^{\omega}$ is 
${\sigma}$-second countable if and only if  $\nw(X)={\omega}$. 
The next theorem shows that $MA(Cohen)$ 
is not enough to get this equivalence.
Note that for a first countable space $X$ we have 
$\w(X)\le |X|$.

\begin{theorem}\label{tm:consistent}
If ZFC is consistent, then so is ZFC + 
MA(Cohen) +``{\em there is a first countable 
0-dimensional Hausdorff space $Y$ such that  $\nw(Y)={\omega}$ and  
 $|Y|<2^{\omega}$ but $Y$ is not
${\sigma}$-second countable
}''.
\end{theorem}

The proof is based on theorem \ref{lm:somewhere-dense} above and theorem 
\ref{tm:second-preserve} below. 
\newcommand{\bdot}{{\dot B}}

\begin{definition}
Given a topological space $Z$, let $\ical(Z)$ be the 
$\sigma$-ideal generated by the nowhere dense subsets of $Z$. 
The elements of $\ical(Z)$ are called {\em first category in $Z$}. 
\end{definition}

The next result is due to Saharon Shelah \cite{Sh}
and it is included here with his kind permission.
\begin{theorem}\label{tm:second-preserve}
If $Z$ is a topological space and  the forcing notion 
$P=\fn(\kappa,2,\omega)$ adds $\kappa$ Cohen
reals to the ground model, then the ideal
 $\ical^{V^P}(Z)$ is generated by $\ical^V(Z)$, that is, 
for each $T\in \ical^{V^P}(Z)$ there is  $T'\in\ical^V(Z)$ with $T\subs T'$.
\end{theorem}
\proof
Since %
$\ical^{V^P}(Z)$ is $\sigma$-generated by the nowhere dense subsets
of $Z$ in $V^P$,
we can assume that $T$ is nowhere dense. Let $\dot{T}$ be a $P$-name of $T$
such that $1_P\force$``{\em $\dot{T}$  is nowhere dense}''.

For each $m\in {\omega}$ define, in $V$, the subset $B_m$
of $Z$ as follows:
\begin{equation}
\tag{$*$}B_m=\{x\in Z:\exists p\in P\ |p|=m\land p\force x\in\dot{T}\}.
\end{equation}
Clearly $1_P\force \dot{T}\subs\bigcup\limits_{m\in {\omega}}B_m $,
so it is enough to show that every $B_m$ is nowhere dense.
Assume on the contrary that there are an open set $U\subs Z$
and $m\in {\omega}$ such that $B_m$ is dense in $U$.

Now, by finite induction, we can define open sets 
$U\supset U_0\supset U_1\dots\supset U_m$ and conditions $q_0,\dots,q_m$
with pairwise disjoint domains such that for each $j\le m$ and for each
$f\in P$ if $\dom(f)=\bigcup\limits_{i<j}\dom(q_j)$ then
\begin{equation}
\tag{\dag}f\cup q_j\force \dot{T}\cap U_j=\empt.
\end{equation}

Since $B_m$ is dense in $U$ and $U_m\subs U$ there is
$x\in B_m\cap U_m$. Then, by the definition of $B_m$, we have
a condition $p\in P$ with $|p|=m$ such that 
$p\force x\in \dot{T}$. But the domains of the $q_j$ are pairwise
disjoint, so there is $j\le m$ with $\dom(p)\cap \dom(q_j)=\empt$.
Thus $p\cup q_j\in P$. Let $f=p\rest \bigcup\limits_{i<j}\dom(q_i)$.
By (\dag) $f\cup q_j\force \dot{T}\cap U_j=\empt$,
so $p\cup q_j\force \dot{T}\cap U_m=\empt$ as well.
But this contradicts $x\in U_m$ and $p\force x\in \dot{T}$, and thus
 the theorem is proved.	
\eproof

\Proof{theorem \ref{tm:consistent}}
Using  again \cite[theorem 3.1]{JSSz}  we have a  
$0$-dimensional Hausdorff space $X$ such that 
${\chi}(X)\nw(X)={\omega}$, but $\w(X)=2^{\omega}$. Let $Y=X^{\omega}$
and $Z=\D(X)^\omega$.
Then add $\kappa>\cont$ Cohen reals to the ground model.
In the generic extension clearly ${\chi}(Y)=\nw(Y)={\omega}$
will remain valid.
 
Since $Z\notin\ical^V(Z)$ by the Baire Category theorem, we have
$Z\notin \ical^{V^P}(Z)$ as well by theorem \ref{tm:second-preserve}.  

Therefore, if $Y=\bigcup\{Y_k:k<\omega\}$ holds in $V^P$,
then $Y_k\notin \ical^{V^P}(Z)$ for some $k\in\omega$, i.e.. 
 $Y_k$ is somewhere dense in $Z$. 
Thus there is a natural number $n\in {\omega}$ 
and a function $f:n\to X$ such that   
\begin{equation}
\tag{\dag}\forall f'\in\Fin({\omega}\setm n,X)
\ \exists g\in Y_k\ f\cup f'\subs g.
\end{equation}

But $(\dag)$ implies that $Y_k$ is  also dense in
the basic open set $(V_f)^{V^P}=\{g\in X^{\omega}\cap V^P:f\subs g\}$ of
$(\D(X)^\omega)^{V^P}$. Thus, applying lemma \ref{lm:somewhere-dense}
in $V^P$ we have $\w(Y_k)=\w(X)>\omega$.
Thus, in $V^P$, $\operatorname{MA}(Cohen)$ holds, and still the 
first countable, 0-dimensional $T_2$ space $Y$ is not
$\sigma$-second countable, though $\w(Y)=(\cont)^V<\kappa=(\cont)^{V^P}$.
\Eproof

\section{Examples for higher cardinals}

In this section we generalize the constructions of \cite[theorem 3.1]{JSSz}
and \ref{tm:need-forcing} for cardinals greater than ${\omega}$.

\begin{theorem}
\label{tm:gen-con}
Let ${\kappa}$ and ${\lambda}$ be  cardinals, 
$\cf({\kappa})={\kappa}$. Then there is a 
0-dimensional Hausdorff space
$\dkl$ such that ${\chi}(\dkl)={\kappa}$, 
$\nw(\dkl)=\lvk$ and $\w(\dkl)=\lk$.
\end{theorem}

\proof
For each $f\in \fn({\kappa},{\lambda},{\kappa})$ put
$U_f=\{g\in \kl:f\subs g\}$.
Write $U(g,{\alpha})=U_{g\rest {\alpha}}$ for 
$g\in \kl$ and ${\alpha}<{\kappa}$.
For $g\ne h\in \kl$ define 
${\Delta}(g,h)=\min\{{\alpha}:g({\alpha})\ne h({\alpha})\}$.
Consider the topological space $\cckl=\<\kl,{\tau}\>$ that has 
as a  base  $\{U_f:f\in \vkl\}$.
Let
\begin{gather}\notag
Y=\{g\in \kl:\exists {\alpha}<{\kappa}\ ( 
\ g({\beta})=0\text{ iff }{\ \beta} \ge {\alpha}\  )\}\\
\intertext{and}\notag
Z=\{g\in \kl:0\notin\ran(g)\}.
\end{gather} 

Clearly $Y$ and $Z$ are disjoint, 
$|Y|=\lvk$, $|Z|=\lk$, $Z$ is closed and nowhere dense in 
$\cckl$.
Let $X=Y\cup Z$. Our required space will be
$\dkl=\<X,{\rho}\>$, where ${\rho}$ refines the topology ${\tau}_X$.
To define ${\rho}$  put
$$X_g=\bigcup\{U_f:\exists {\alpha}<{\kappa}\ 
f=g\rest {\alpha}\cup\{\<{\alpha},0\>\}\}$$
for $g\in Z$ and  $X_g=\empt$ for $g\in Y$.
For $g\in X$ and ${\alpha}<{\kappa}$ let
  $V(g,{\alpha})=(U(g,{\alpha})\setm  X_g)\cap X$.
Let the neighbourhood base of $g\in X$ in ${\rho}$ be 
$$
\bcal_g=\{V(g,{\alpha}):{\alpha}<{\kappa}\}.
$$

First we note that   
$\bcal=\bigcup\{\bcal_x:x\in X\}$ is a base of a topology
because
\begin{equation}
\tag{\dag}\forall {\alpha}<{\kappa}\ \forall h\in V(g,{\alpha})\setm \{g\}
\  U(h,{\Delta}(g,h)+1)\cap X\subs V(g,{\alpha}).
\end{equation}

Since $V(g,{\alpha})\cap Y=U(g,{\alpha})\cap Y$ and 
$V(h,{\alpha})\cap Z=U(h,{\alpha})\cap Z$ for each 
$g\in Y$,  $h\in Z$ and ${\alpha}<{\kappa}$
we have that $\<Y,{\tau}\>=\<Y,{\rho}\>$ and $\<Z,{\tau}\>=\<Z,{\rho}\>$.
Thus 
$$
\nw(\<X,{\rho}\>)=\nw(\<Y,{\rho}\>)+\nw(\<Z,{\rho}\>)\le
\w(\<Y,{\tau}\>)+\w(\<Z,{\tau}\>)\le \w(\<\cckl\>)= \lvk.
$$
Obviously ${\chi}(\<X,{\rho}\>)={\kappa}$.

Finally we show that $\w(\<X,{\rho}\>)=\lk$.
This will follow if we show that 
$\bcal$ is an irreducible base for $X$, by lemma
\cite[2.6]{JSSz}.
We claim that $\{\bcal_x:x\in X\}$ is an irreducible decomposition
of the base $\bcal$ 
(see definition \cite[2.3]{JSSz}). 
Since  $Y$ is discrete in ${\rho}$ 
 it is enough to show that if 
$g\ne h$ are from $Z$ with  $g\in V(h,{\alpha})$ for some
 ${\alpha}<{\kappa}$ then $V(h,{\alpha})\not\subs V(g,0)$. 
Let ${\delta}={\Delta}(g,h)$.
Then 
$U(g,{\delta}+1)\subs V(h,{\alpha})$ by $(\dag)$.
Consider the   element $y$ of $Y$ defined by the formulas 
$y\rest {\delta}+1=g\rest {\delta}+1$ and $y({\delta}+1)=0$.
Then $y\in X_g$ and  so $y\notin V(g,0)$. 
On the other hand $y\in U(g,{\delta}+1)$, so $V(h,{\alpha})\not\subs V(g,0)$.
\eproof

\begin{lemma}\label{lm:bigsomewhere-dense}
If  $X$ is any topological space and $Z\subs X^{\kappa}$,  
${\alpha}<{\kappa}$,  $f:{\alpha}\to X$ are such that 
\begin{equation}
\tag{\ddag}\forall f'\in\Fin({\kappa}\setm {\alpha},X)
\ \exists g\in Z\ f\cup f'\subs g,
\end{equation}
then $\w(Z)\ge\w(X)$.
\end{lemma}

The proof is similar to  that of lemma \ref{lm:somewhere-dense}, so
we omit it.

Let us recall that given a cardinal ${\mu}$  the {\em 
Singular Cardinal Hypothesis $($SCH$)$} is said 
to hold below ${\mu}$ provided
${\nu}^{\cf({\nu})}=2^{\cf({\nu})}{\nu}^+$
for each singular cardinal ${\nu}\le {\mu}$.
By \cite[lemma 1.8.1]{Jech}, if $\mu$ is regular and  SCH holds below $\mu$, 
then $\log(\mu^+)=\min\{\nu:2^\nu\ge\mu^+\}$ is also regular.
It is well-known that the failure of SCH requires the consistency 
of large cardinals, therefore the assumption of our next lemma is quite 
reasonable. Also note that $SCH$ trivially holds below $\aleph_{\omega}$.

\begin{theorem}
\label{tm:sch}
If ${\mu}$ and $\log({\mu^+})$ are both regular cardinals then
there is a 0-dimensional $T_2$  space $Y$ such that 
${\chi}(Y)\nw(Y)\le{\mu}$, but $Y$ is not the union of ${\mu}$ subspaces
of weight ${\mu}$.
\end{theorem}

\proof
Let ${\rho}=\log(\mu^+)$.
Applying  theorem \ref{tm:gen-con} for ${\kappa}={\rho}$ and ${\lambda}=2$
we get a space $X$ with 
${\chi}(X)\nw(X)= 2^{<{\rho}}\le {\mu}<\w(X)$. 
Let $Y=X^{\rho}$ and consider a partition 
$Y=\bigcup_{{\alpha}<{\rho}} Y_{\alpha}$. 
Since $\kappa={\rho}$ is regular, applying the technique of the 
standard proof of the  Baire Category theorem we can see  that  
the space $C^\rho_{|X|}$ is not the union of 
$\rho$ nowhere dense subspaces. 
Therefore there  are  ordinals ${\alpha},{\xi}<{\rho}$ and a  
function $f:{\alpha}\to X$ such that 
\begin{equation}
\tag{\ddag'}\forall f'\in\fn({\rho}\setm {\alpha},X,\rho)
\ \exists g\in Y_{\xi}\ f\cup f'\subs g.
\end{equation}
But  (\ddag') clearly implies (\ddag). 
So, by lemma \ref{lm:bigsomewhere-dense},
$\w(Y_{\xi})\ge \w(X)>{\mu}$. Thus $Y$ satisfies the requirements.
The theorem is proved.
\eproof

\begin{problem} $($ZFC$)$
If $\mu$ is a cardinal, is there a  space 
$X$ such that $\chi(X)\nw(X)=\mu$
(or just $\nw(X)={\mu}$), but $X$ is
not the union of $<2^\mu$ many subspaces of
weight $<2^\mu$ $($or just $\w(X)=2^\mu)$? 
\end{problem}

\begin{remark}
The simplest case of Problem 1 left open by theorem \ref{tm:gen-con} 
is that $2^{\omega}={\omega}_2$, $2^\oo={\omega}_3$ and
${\mu}=\oo$. Indeed, if $\dkl$ is the space constructed 
in theorem \ref{tm:gen-con} and  $\nw(\dkl)=\lvk\le\oo$, 
then ${\lambda}\le\oo$ and ${\kappa}\le{\omega}$. So 
$\w(\dkl)=\lk\le \oo^{\omega}={\omega}_2<2^{\oo}$. 
\end{remark}

\end{document}